**[1,2]Манарбек М.,**
PhD студент, ORCID ID: 0009-0006-6879-8356,
e-mail: manarbek@math.kz
**[2]Тлеуханова Н. Т.,**
профессор, физика-математика ғылымдарының докторы,
ORCID ID: 0000-0002-4133-7780,
e-mail: tleukhanova@rambler.ru
**[2]Мусабаева Г. К.,**
PhD, ORCID ID: 0000-0003-2368-8955,
e-mail: musabaevaguliya@mail.ru

[1]Математика және математикалық модельдеу институты, Алматы қ., Қазақстан
[2]Л.Н. Гумилев атындағы Еуразия ұлттық университеті, Астана қ., Қазақстан

# АНИЗОТРОПТЫ ГРАНД ЛОРЕНЦ КЕҢІСТІКТЕРІ ЖӘНЕ ОЛАРДЫҢ ҚАСИЕТТЕРІ

**Аңдатпа**

Бұл мақалада жаңа анизотропты гранд Лоренц кеңістіктері анықталып, олардың қасиеттері зерттеледі. Бұл кеңістіктер әртүрлі функционалдық кеңістіктерді зерттеу үшін біртұтас параметрді қамтамасыз ететін жаңа құрылым болып табылады. Гранд кеңістіктерді қарастыру әсіресе параметрлердің шектік жағдайларын зерттеу үшін маңызды және осы мәселеде жаңа нәтижелерге қол жеткізуге болады. Классикалық кеңістіктерде шектік параметрлерді зерттеу әр уақытта мүмкін бола бермейді. Соңғы жылдары функционалдық кеңістіктер мәселелерінде гранд Лебег кеңістіктері және олардың жалпылаулары кеңінен зерттеліп жүр. Бұл кеңістіктер классикалық Лоренц және гранд Лоренц кеңістіктерінің жалпыламалары болады. Мақалада гранд анизотропты Лоренц кеңістіктерінің анықтамасы, осы кеңістіктегі негізгі бағалаулар көрсетілді, енгізу теоремалары дәлелденді, параметрлер бойынша енгізу шарттары шығарылды. Алынған нәтижелер теориялық тұрғыда ғана емес, сонымен қатар қолданбалы есептерде де маңызды рөл атқара алады.

**Тірек сөздер:** Лоренц кеңістіктері, гранд Лоренц кеңістіктері, енгізу теоремалары, теңсіздіктер, анизотропты кеңістіктер.

## Кіріспе

Соңғы жиырма жыл ішінде гранд Лебег кеңістіктер теориясы заманауи талдаудың интенсивті дамып келе жатқан бағыттарының бірі. Бұл кеңістіктер $L^{p)}(\Omega)$ деп белгіленіп, $\Omega$ облысында анықталған және келесі шарттарды қанағаттандыратын барлық өлшенетін $f$ функциялар жиыны:

$$\|f\|_{L^{p),\theta}(\Omega)} := \sup_{0<\varepsilon<p-1} \varepsilon^{\theta} \|f\|_{L^{p-\varepsilon}(\Omega)} < \infty, |\Omega| < \infty$$

Гранд Лебег кеңістіктері 1992 жылы Иваньек және Сбордонмен енгізілген болатын [1]. Олар бұл кеңістіктерді Якобианның интегралдануы мәселесін минималды болжамдармен шешу үшін қолданды.

Кейінірек, $|\Omega|=1$ болған жағдайда гранд Лебег кеңістіктерінің альтернативті сипаттамасы [2]-[3], жұмыстарында берілген:

$$\|f\|_{L^{p,\theta}(\Omega)} \approx \sup_{0<s<1}(1-\ln s)^{-\frac{\theta}{p}}\left(\int_s^1 (f^*(t))^p \frac{dt}{t}\right)^{\frac{1}{p}}.$$

Гранд Лебег кеңістіктері дифференциалдық теңдеулерге байланысты кейбір мәселелерде [4]-[7] қолданылған. Сонымен қатар, бұл кеңістіктер сызықтық емес теңдеулерді шешу үшін қолайлы орта ретінде танылды. Гранд Лебег кеңістіктеріндегі операторлар теориясы соңғы жылдары кеңінен зерттелді. Бұл саладағы көптеген нәтижелер осы тақырыпқа арналған В. Кокилашвили, А. Месхи, У. Рафейро авторларының [8]-[9] монографияларында жинақталған. $Lpq$ Лоренц кеңістіктерін 1950 жылдары Г. Лоренц [10]-[11] енгізген, және олардың маңыздылығы, интерполяция теориясында Лебег кеңістіктерін толықтыратындығында.

Классикалық деректер үшін С. Беннетт и Р. С. Шарплидің [12] жұмыстарына сілтеме жасауға болады, ал заманауи деректер үшін Р. Е. Кастильо и Х. С. Чапарро [13] еңбектерін қарастыруға болады. Сонымен қатар классикалық кеңістіктер үшін Е.Д. Нұрсұлтанов [14]-[17] жұмыстарда одан әрі дамытты. Анизотропты және аралас метрикалы Лоренц кеңістіктеріндегі теңсіздіктер жайлы Н. Т. Тлеуханова және Г. К. Мусабаева [18]-[19] мақалаларында зерттелді. Е.Д. Нұрсұлтанов, У. Рафейро, және Д. Сураган [20] жұмысында гранд Лоренц кеңістіктерін анықтады. Лоренц–Карамата кеңістігі мен гранд Лоренц кеңістігінің басқа нұсқасымен салыстырулар жүргізілді.

Бұл мақалада Юнг пен О'Нэйл теңсіздіктері енгізілген гранд Лоренц кеңістіктерінде дәлелденді, бұл теңсіздіктер Харди-Литлвуд-Соболев түріндегі теңсіздікті шығаруға мүмкіндік береді. [20] мақаласында гранд Лебег, анизотропты емес гранд Лоренц кеңістіктерінің әр түрлі нұсқалары анықталып, қасиеттері зерттелген. Ал берілген мақалада анизотропты гранд Лоренц кеңістігі алғаш енгізіліп тұр. Жаңадан енгізілген гранд Лоренц кеңістіктерінің маңызды қасиеті – классикалық Лоренц кеңістіктері оң және теріс $\theta$ мәндері бар гранд Лоренц кеңістіктерінің арасында орналасатынында:

$$GL_{p,q}^{-\theta}(\Omega) \hookrightarrow L_{p,q}(\Omega) \hookrightarrow GL_{p,q}^{\theta}(\Omega)$$

Гранд Лоренц кеңістіктерінде интерполяция мәселелеріне тоқталды. Сондай-ақ гранд Лоренц кеңістіктері үшін Кёте дуалдылығын талқылап, одан гранд Лебег кеңістіктеріндегі Кёте дуалды кеңістік теоремасы алынды.

**Материалдар мен әдістер**

Жұмыста екі өлшемді $GL_{\bar{p},\bar{q}}^{\bar{\theta}}(\Omega)$ анизотропты гранд Лоренц кеңістіктері, олардың қасиеттері қарастырылады. Мұнда $\bar{p} = (p_1, p_2)$, $0 < p_1, p_2 \leq \infty$,

$\frac{1}{p_1} + \frac{1}{p_1'} = 1, \frac{1}{p_2} + \frac{1}{p_2'} = 1$ болған жағдайларды $0 < \bar{p} \leq \infty, \frac{1}{\bar{p}} + \frac{1}{\bar{p}'} = 1$ деп жазамыз.

**Анықтама 1.** $\bar{\theta} = (\theta_1, \theta_2)$, $-\infty < \bar{\theta} \leq \infty$, $\bar{p} = (p_1, p_2)$, $0 < \bar{p} \leq \infty$, $0 < \bar{q} \leq \infty$, $\bar{q} = (q_1, q_2)$ және $\bar{\theta} \in \mathbb{R}^2, |\Omega| = 1$ болсын. $GL_{\bar{p},\bar{q}}^{\bar{\theta}}(\Omega)$ анизотропты гранд Лоренц кеңістіктері деп келесі квазинормалары ақырлы болатын $f$ өлшенетін функциялар жиынын айтамыз:

$0 < q < \infty$ кезінде:

$$\|f\|_{GL_{\bar{p},\bar{q}}^{\bar{\theta}}(\Omega)} = \begin{cases} \theta_1 \geq 0, \theta_2 \geq 0 \text{ үшін:} \\ \sup_{0<\varepsilon_1,\varepsilon_2 \leq 1} \varepsilon_1^{\theta_1} \varepsilon_2^{\theta_2} \left( \int_0^1 \left( \int_0^1 \left( t_1^{\frac{1}{p_1}+\varepsilon_1} t_2^{\frac{1}{p_2}+\varepsilon_2} f^{*_1,*_2}(t_1,t_2) \right)^{q_1} \frac{dt_1}{t_1} \right)^{\frac{q_2}{q_1}} \frac{dt_2}{t_2} \right)^{\frac{1}{q_2}}, \\ \theta_1 > 0, \theta_2 > 0, \quad _1 = \infty, p_2 = \infty, \text{ үшін:} \\ \sup_{0<\varepsilon_1,\varepsilon_2 \leq 1} \varepsilon_1^{\theta_1} \varepsilon_2^{\theta_2} \left( \int_0^1 \left( \int_0^1 \left( t_1^{\varepsilon_1} t_2^{\varepsilon_2} f^{*_1,*_2}(t_1,t_2) \right)^{q_1} \frac{dt_1}{t_1} \right)^{\frac{q_2}{q_1}} \frac{dt_2}{t_2} \right)^{\frac{1}{q_2}}, \\ \theta_1 < 0, \theta_2 < 0 \text{ үшін:} \\ \inf_{\substack{0<\varepsilon_1 \leq \frac{1}{p_1} \\ 0<\varepsilon_2 \leq \frac{1}{p_2}}} \varepsilon_1^{\theta_1} \varepsilon_2^{\theta_2} \left( \int_0^1 \left( \int_0^1 \left( t_1^{\frac{1}{p_1}-\varepsilon_1} t_2^{\frac{1}{p_2}-\varepsilon_2} f^{*_1,*_2}(t_1,t_2) \right)^{q_1} \frac{dt_1}{t_1} \right)^{\frac{q_2}{q_1}} \frac{dt_2}{t_2} \right)^{\frac{1}{q_2}}, \end{cases} \quad (1)$$

және $q = \infty$ кезінде:

$$\|f\|_{GL_{\bar{p},\infty}^{\bar{\theta}}(\Omega)} = \begin{cases} \theta_1 \geq 0, \theta_2 \geq 0 \text{ үшін:} \\ \sup_{0<\varepsilon_1,\varepsilon_2<1} \sup_{0<t_1,t_2<1} \varepsilon_1^{\theta_1} \varepsilon_2^{\theta_2} t_1^{\frac{1}{p_1}+\varepsilon_1} t_2^{\frac{1}{p_2}+\varepsilon_2} f^{*_1,*_2}(t_1,t_2), \\ \theta_1 > 0, \theta_2 > 0, p_1 = \infty, p_2 = \infty \text{ үшін:} \\ \sup_{0<\varepsilon_1,\varepsilon_2<1} \sup_{0<t_1,t_2<1} \varepsilon_1^{\theta_1} \varepsilon_2^{\theta_2} t_1^{\varepsilon_1} t_2^{\varepsilon_2} f^{*_1,*_2}(t_1,t_2), \\ \theta_1 < 0, \theta_2 < 0 \text{ үшін:} \\ \inf_{0<\varepsilon_1,\varepsilon_2<\frac{1}{p_1},\frac{1}{p_2}} \sup_{0<t_1,t_2<1} \varepsilon_1^{\theta_1} \varepsilon_2^{\theta_2} t_1^{\frac{1}{p_1}-\varepsilon_1} t_2^{\frac{1}{p_2}-\varepsilon_2} f^{*_1,*_2}(t_1,t_2) \end{cases} \quad (2)$$

**Ескерту.** $GL_{\bar{p},\bar{q}}^{\bar{\theta}}(\Omega)$ Гранд анизотропты Лоренц кеңістіктерінің классикалық анизотропты $L_{\bar{p},\bar{q}}(\Omega)$ Лоренц кеңістіктерінен ерекшелігі $GL_{\bar{p},\infty}^{\bar{\theta}}(\Omega)$ кеңістіктер шкаласын қарастыруға мүмкіндік береді. Егер $\varepsilon_1=\varepsilon_2=0$ болса, біз классикалық Лоренц кеңістігін аламыз. Қалған жағдайларда бұл жаңа кеңістік болып саналады.

**Нәтижелер мен талқылау**

**Теорема 1.** $\bar{\theta} > 0$ болған кезде келесі енгізулер орындалады:

$$GL_{\bar{p},\bar{q}}^{-\bar{\theta}}(\Omega) \hookrightarrow L_{\bar{p},\bar{q}}(\Omega) \hookrightarrow GL_{\bar{p},\bar{q}}^{\bar{\theta}}(\Omega). \qquad (3)$$

**Дәлелдеу.** $\theta_1, \theta_2 > 0$ болған жағдайда $L_{\bar{p},\bar{q}}(\Omega) \hookrightarrow GL_{\bar{p},\bar{q}}^{\bar{\theta}}(\Omega)$ енгізуі орындалатынын көрсетейік:

$0 < \varepsilon_i < 1,\ i \in (1,2)$ және $t^{\frac{1}{\bar{p}}} + \bar{\varepsilon} < t^{\frac{1}{\bar{p}}},\ 0 < t < 1$ екенін ескере отырып, (1) анықтамасын қолданып, төмендегі теңсіздіктерді аламыз:

$$\|f\|_{GL_{\bar{p},\bar{q}}^{\bar{\theta}}(\Omega)} = \sup_{0<\varepsilon_1,\varepsilon_2\leq 1} \varepsilon_1^{\theta_1}\varepsilon_2^{\theta_2} \left( \int_0^1 \left( \int_0^1 \left( t_1^{\frac{1}{p_1}+\varepsilon_1} t_2^{\frac{1}{p_2}+\varepsilon_2} f^{*_1,*_2}(t_1,t_2) \right)^{q_1} \frac{dt_1}{t_1} \right)^{\frac{q_2}{q_1}} \frac{dt_2}{t_2} \right)^{\frac{1}{q_2}}$$

$$\leq C \left( \int_0^1 \left( \int_0^1 \left( t_2^{\frac{1}{p_2}} t_1^{\frac{1}{p_1}} f^{*_1*_2}(t_1,t_2) \right)^{q_1} \frac{dt_1}{t_1} \right)^{\frac{q_2}{q_1}} \frac{dt_2}{t_2} \right)^{\frac{1}{q_2}} = C\|f\|_{L_{\bar{p},\bar{q}}(\Omega)}.$$

Енді $GL_{\bar{p},\bar{q}}^{-\bar{\theta}}(\Omega) \hookrightarrow L_{\bar{p},\bar{q}}(\Omega)$ енгізуін дәлелдейік. $0 < t_1,t_2 \leqslant 1, t^{-\varepsilon} > 1$ екенін ескеріп және инфинумның қасиеттерін қолдана отырып, төмендегі теңсіздікті аламыз:

$$\|f\|_{L_{\bar{p},\bar{q}}(\Omega)} = \left( \int_0^1 \left( \int_0^1 \left( t_2^{\frac{1}{p_2}} t_1^{\frac{1}{p_1}} f^{*_1*_2}(t_1,t_2) \right)^{q_1} \frac{dt_1}{t_1} \right)^{\frac{q_2}{q_1}} \frac{dt_2}{t_2} \right)^{\frac{1}{q_2}}$$

$$\leq \left( \int_0^1 \left( \int_0^1 \left( t_2^{\frac{1}{p_2}-\varepsilon_2} t_1^{\frac{1}{p_1}-\varepsilon_1} f^{*_1*_2}(t_1,t_2) \right)^{q_1} \frac{dt_1}{t_1} \right)^{\frac{q_2}{q_1}} \frac{dt_2}{t_2} \right)^{\frac{1}{q_2}},$$

$$\|f\|_{L_{\bar{p},\bar{q}}(\Omega)} \leq \inf_{\substack{0<\varepsilon_1 \leq \frac{1}{p_1} \\ 0<\varepsilon_2 \leq \frac{1}{p_2}}} \varepsilon_1^{\theta_1} \varepsilon_2^{\theta_2} \left( \int_0^1 \left( \int_0^1 \left( t_1^{\frac{1}{p_1}-\varepsilon_1} t_2^{\frac{1}{p_2}-\varepsilon_2} f^{*_1,*_2}(t_1,t_2|\Omega|) \right)^{q_1} \frac{dt_1}{t_1} \right)^{\frac{q_2}{q_1}} \frac{dt_2}{t_2} \right)^{\frac{1}{q_2}}$$

$$\leq C \|f\|_{GL_{\bar{p},\bar{q}}^{\bar{\theta}}(\Omega)}$$

$\bar{p} = \bar{q}$ дербес жағдайында (3) енгізулерінен келесі енгізулер орынды екенін көреміз:

$$GL_{\bar{p},\bar{p}}^{-\bar{\theta}}(\Omega) \hookrightarrow L_{\bar{p},\bar{p}}(\Omega) \hookrightarrow GL_{\bar{p},\bar{p}}^{\bar{\theta}}(\Omega).$$

**Теорема 2.** $0 < p_1, p_2 < \infty$ болсын. Онда $\theta_1, \theta_2 > 0$ үшін келесі бағалаулар орындалады:

1. $\|f\|_{GL_{\bar{p},\bar{q}}^{\bar{\theta}}(\Omega)} \asymp \sup_{t_1,t_2>0} \frac{t_1^{\frac{1}{p_1}} t_2^{\frac{1}{p_2}}}{|\ln t_1|^{\theta_1}|\ln t_2|^{\theta_2}} f^{*_1,*_2}(t_1,t_2),$

2. $\|f\|_{GL_{\bar{p},\bar{q}}^{\bar{\theta}}(\Omega)} \lesssim \left( \int_0^1 \left( \int_0^1 \left( \frac{t_1^{\frac{1}{p_1}}, t_2^{\frac{1}{p_2}}}{|\ln t_1|^{\theta_1}|\ln t_2|^{\theta_2}} f^{*_1,*_2}(t_1,t_2) \right)^{q_1} \frac{dt_1}{t_1} \right)^{\frac{q_2}{q_1}} \frac{dt_2}{t_2} \right)^{\frac{1}{q_2}},$

3. $\|f\|_{GL_{\bar{p},\bar{q}}^{-\bar{\theta}}(\Omega)} \gtrsim \left( \int_0^1 \left( \int_0^1 \left( t_1^{\frac{1}{p_1}}|\ln t_1|^{\theta_1} t_2^{\frac{1}{p_2}}|\ln t_2|^{\theta_2} f^{*_1,*_2}(t_1,t_2) \right)^{q_1} \frac{dt_1}{t_1} \right)^{\frac{q_2}{q_1}} \frac{dt_2}{t_2} \right)^{\frac{1}{q_2}}.$

**Дәлелдеу.** Келесі функцияны қарастырырып, экстремумге зерттейік:
$\varphi(\bar{\varepsilon}) = \left( \varepsilon_1^{\theta_1} t_1^{\frac{1}{p_1}+\varepsilon_1 \text{sign}\theta_1} \right) \left( \varepsilon_2^{\theta_2} t_2^{\frac{1}{p_2}+\varepsilon_2 \text{sign}\theta_2} \right)$, мұндағы $0 < t_1, t_2 \leq 1, 0 < \varepsilon_1, \varepsilon_2, \leq 1$ болсын. Келесі теңдіктерден:

$\frac{\partial \varphi}{\partial \varepsilon_1} = \left( \theta_1 \varepsilon_1^{\theta_1-1} t_1^{\frac{1}{p_1}+\varepsilon_1 \text{sign}\theta_1} + \varepsilon_1^{\theta_1} t_1^{\frac{1}{p_1}+\varepsilon_1 \text{sign}\theta_1} \text{sign}\theta_1 \ln t_1 \right) \left( \varepsilon_2^{\theta_2} t_2^{\frac{1}{p_2}+\varepsilon_2 \text{sign}\theta_2} \right) = 0,$

$\frac{\partial \varphi}{\partial \varepsilon_2} = \left( \theta_2 \varepsilon_2^{\theta_2-1} t_2^{\frac{1}{p_2}+\varepsilon_2 \text{sign}\theta_2} + \varepsilon_2^{\theta_2} t_2^{\frac{1}{p_2}+\varepsilon_2 \text{sign}\theta_2} \text{sign}\theta_2 \ln t_2 \right) \left( \varepsilon_1^{\theta_1} t_1^{\frac{1}{p_1}+\varepsilon_1 \text{sign}\theta_1} \right) = 0.$

$\varepsilon_1 = \frac{|\theta_1|}{|\ln t_1|}$ және $\varepsilon_2 = \frac{|\theta_2|}{|\ln t_2|}$ болғанда $\varphi(\bar{\varepsilon})$ функциясының супремумы бар. Демек $\theta_1, \theta_2 > 0$, $\bar{t}^{\bar{\varepsilon}} \leq 1$, $0 < t_1, t_2 \leq 1$ үшін:

$$\sup_{0<\bar{\varepsilon}\leq 1} \varphi(\bar{\varepsilon}) = \sup_{0<\bar{\varepsilon}\leq 1} \left( \varepsilon_1^{\theta_1} t_1^{\frac{1}{p_1}+\varepsilon_1 sign\theta_1} \right) \left( \varepsilon_2^{\theta_2} t_2^{\frac{1}{p_2}+\varepsilon_2 sign\theta_2} \right) \asymp \frac{t_1^{\frac{1}{p_1}} t_2^{\frac{1}{p_2}}}{|ln\ t_1|^{\theta_1}|ln\ t_2|^{\theta_2}},$$

ал $\theta_1, \theta_2 < 0$ үшін

$$\inf_{0<\bar{\varepsilon}\leq 1} \varphi(\bar{\varepsilon}) = \inf_{0<\bar{\varepsilon}\leq 1}\left(\varepsilon_1^{\theta_1} t_1^{\frac{1}{p_1}+\varepsilon_1 \text{sign}\theta_1}\right)\left(\varepsilon_2^{\theta_2} t_2^{\frac{1}{p_2}+\varepsilon_2 \text{sign}\theta_2}\right) \asymp \frac{t_1^{\frac{1}{p_1}} t_2^{\frac{1}{p_2}}}{|\ln t_1|^{\theta_1}|\ln t_2|^{\theta_2}}.$$

Бұл дегеніміз

$$\|f\|_{GL_{\bar{p},\infty}^{\bar{\theta}}(\Omega)} = \sup_{0<\varepsilon_1,\varepsilon_2<1} \sup_{t_1,t_2>0} \varepsilon_1^{\theta_1} t_1^{\frac{1}{p_1}+\varepsilon_1} \varepsilon_2^{\theta_2} t_2^{\frac{1}{p_2}+\varepsilon_2} f^{*_1,*_2}(t_1,t_2)$$

$$= \sup_{t_1,t_2>0}\left(\sup_{0<\varepsilon_1,\varepsilon_2\leq 1}\varepsilon_1^{\theta_1} t_1^{\frac{1}{p_1}+\varepsilon_1} \varepsilon_2^{\theta_2} t_2^{\frac{1}{p_2}+\varepsilon_2}\right) f^{*_1,*_2}(t_1,t_2)$$

$$\asymp \sup_{t_1,t_2>0} \frac{t_1^{\frac{1}{p_1}} t_2^{\frac{1}{p_2}}}{|\ln t_1|^{\theta_1}|\ln t_2|^{\theta_2}} f^{*_1,*_2}(t_1,t_2).$$

Енді екінші бағалауды дәлелдейік:

$$\|f\|_{GL_{\bar{p},\bar{\tau}}^{\bar{\theta}}(\Omega)} = \sup_{0<\varepsilon_1,\varepsilon_2\leq 1} \varepsilon_1^{\theta_1}\varepsilon_2^{\theta_2}\left(\int_0^1\left(\int_0^1\left(t_1^{\frac{1}{p_1}+\varepsilon_1} t_2^{\frac{1}{p_2}+\varepsilon_2} f^{*_1,*_2}(t_1,t_2)\right)^{\tau_1}\frac{dt_1}{t_1}\right)^{\frac{\tau_2}{\tau_1}}\frac{dt_2}{t_2}\right)^{\frac{1}{\tau_2}}$$

$$\leq \left(\int_0^1\left(\int_0^1 \sup_{0<\varepsilon_1,\varepsilon_2\leq 1}\varepsilon_1^{\theta_1}\varepsilon_2^{\theta_2}\left(t_1^{\frac{1}{p_1}+\varepsilon_1} t_2^{\frac{1}{p_2}+\varepsilon_2} f^{**_1,*_2}(t_1,t_2)\right)^{\tau_1}\frac{dt_1}{t_1}\right)^{\frac{\tau_2}{\tau_1}}\frac{dt_2}{t_2}\right)^{\frac{1}{\tau_2}}$$

$$\asymp \left(\int_0^1\left(\int_0^1\left(\frac{t_1^{\frac{1}{p_1}} t_2^{\frac{1}{p_2}}}{|\ln t_1|^{\theta_1}|\ln t_2|^{\theta_2}} f^{*_1,*_2}(t_1,t_2)\right)^{\tau_1}\frac{dt_1}{t_1}\right)^{\frac{\tau_2}{\tau_1}}\frac{dt_2}{t_2}\right)^{\frac{1}{\tau_2}}.$$

Үшінші бағалауды дәлелдеу үшін 1 пунктті ескереміз:

$$\|f\|_{GL_{\bar{p},\bar{\tau}}^{-\bar{\theta}}(\Omega)} = \inf_{\substack{0<\varepsilon_1<\frac{1}{p_1} \\ 0<\varepsilon_2<\frac{1}{p_2}}} \varepsilon_1^{\theta_1}\varepsilon_2^{\theta_2}\left(\int_0^1\left(\int_0^1\left(t_1^{\frac{1}{p_1}-\varepsilon_1} t_2^{\frac{1}{p_2}-\varepsilon_2} f^{*_1,*_2}(t_1,t_2)\right)^{\tau_1}\frac{dt_1}{t_1}\right)^{\frac{\tau_2}{\tau_1}}\frac{dt_2}{t_2}\right)^{\frac{1}{\tau_2}}$$

$$\geq \left( \int_0^1 \left( \int_0^1 \inf_{\substack{0<\varepsilon_1<\frac{1}{p_1}\\0<\varepsilon_2<\frac{1}{p_2}}} \varepsilon_1^{\theta_1} \varepsilon_2^{\theta_2} \left( t_1^{\frac{1}{p_1}-\varepsilon_1} t_2^{\frac{1}{p_2}-\varepsilon_2} f^{*_1,*_2}(t_1,t_2) \right)^{\tau_1} \frac{dt_1}{t_1} \right)^{\frac{\tau_2}{\tau_1}} \frac{dt_2}{t_2} \right)^{\frac{1}{\tau_2}}$$

$$\asymp \left( \int_0^1 \left( \int_0^1 t_1^{\frac{1}{p_1}} |\ln t_1|^{\theta_1} t_2^{\frac{1}{p_2}} |\ln t_2|^{\theta_2} f^{*_1,*_2}(t_1,t_2) \right)^{\tau_1} \frac{dt_1}{t_1} \right)^{\frac{\tau_2}{\tau_1}} \frac{dt_2}{t_2} \right)^{\frac{1}{\tau_2}}.$$

**Мысал 1.** $\bar{\theta} > 0, \bar{\delta} > 0$, $\Omega \subset R^n$ $|\Omega| = 1$ болсын. Егер $f: \Omega \to R^2$,

$f^{*_1,*_2}(t_1,t_2) \asymp \dfrac{|\ln t_1|^{\theta_1-\frac{1}{r_1}-\delta_1} |\ln t_2|^{\theta_2-\frac{1}{r_2}-\delta_2}}{t_1^{\frac{1}{r_1}} t_2^{\frac{1}{r_2}}}$ орындалса, онда келесі функция

анизотропты гранд Лоренц кеңістігінде жататыны белгілі: $f \in GL_{\bar{p},\bar{r}}^{\bar{\theta}}(\Omega)$.

Келесі теоремаларда гранд анизотропты Лоренц $GL_{\bar{p},\bar{q}}^{\bar{\theta}}(\Omega)$ кеңістіктері үшін әртүрлі параметрлер бойынша енгізу шарттары және де басқа қасиеттері дәлелденеді. Бұл кеңістіктердің қасиеттерін зерттеу арқылы функционалдық кеңістіктердің шектік параметрлерін тереңірек түсінуге, сондай-ақ әртүрлі математикалық есептерді шешуге арналған жаңа әдіс-тәсілдерді әзірлеуге мүмкіндік береді.

**Теорема 3.** Егер $\bar{\theta} \leq \bar{s}$ болса, онда $GL_{\bar{p},\bar{q}}^{\bar{\theta}}(\Omega) \hookrightarrow GL_{\bar{p},\bar{q}}^{\bar{s}}(\Omega)$ енгізуі орындалады.

**Дәлелдеу.** $\theta_1 \leq s_1$, $\theta_2 \leq s_2$, $0 < \varepsilon_1, \varepsilon_2 < 1$, $\varepsilon_1^{s_1} < \varepsilon_1^{\theta_1}$, $\varepsilon_1^{s_2} < \varepsilon_2^{\theta_2}$ екенін ескере отырып, біз келесі теңсіздікті аламыз:

$$\|f\|_{GL_{\bar{p},\bar{q}}^{\bar{s}}(\Omega)} = \sup_{0<\varepsilon_1,\varepsilon_2<1} \varepsilon_1^{s_1} \varepsilon_2^{s_2} \left( \int_0^1 \left( \int_0^1 \left( t_1^{\frac{1}{p_1}+\varepsilon_1} t_2^{\frac{1}{p_2}+\varepsilon_2} f^{*_1,*_2}(t_1,t_2) \right)^{q_1} \frac{dt_1}{t_1} \right)^{\frac{q_2}{q_1}} \frac{dt_2}{t_2} \right)^{\frac{1}{q_2}}$$

$$\leq \sup_{0<\varepsilon_1,\varepsilon_2<1} \varepsilon_1^{\theta_1} \varepsilon_2^{\theta_2} \left( \int_0^1 \left( \int_0^1 \left( t_1^{\frac{1}{p_1}+\varepsilon_1} t_2^{\frac{1}{p_2}+\varepsilon_2} f^{*_1,*_2}(t_1,t_2) \right)^{q_1} \frac{dt_1}{t_1} \right)^{\frac{q_2}{q_1}} \frac{dt_2}{t_2} \right)^{\frac{1}{q_2}}$$

$$= \|f\|_{GL_{\bar{p},\bar{q}}^{\bar{\theta}}(\Omega)}.$$

**Теорема 4.** Егер $\bar{q} < \bar{r}$ болса, онда $GL_{\bar{p},\bar{q}}^{\bar{\theta}}(\Omega) \hookrightarrow GL_{\bar{p},\bar{r}}^{\bar{\theta}}(\Omega)$.

**Дәлелдеу.** $\bar{q} < \bar{r}$ кезінде $L_{p,q}(\Omega) \hookrightarrow L_{p,r}(\Omega)$ енгізуі орындалатынын классикалық Лоренц кеңістігі үшін белгілі енгізу теоремасы [20] мен $GL_{\bar{p},\bar{q}}^{\bar{\theta}}(\Omega)$ анизотропты гранд Лоренц кеңістігінің (1) анықтамасынан шығады.

**Теорема 5.** Егер $0 < \bar{\delta} < 1, \bar{\theta} > 0$ болса, онда

$$\|f\|_{GL_{\bar{p},\bar{q}}^{\bar{\theta}}(\Omega)} \asymp \sup_{\substack{0 < \varepsilon_1 \le \delta_1 \\ 0 < \varepsilon_2 \le \delta_2}} \varepsilon_1^{\theta_1} \varepsilon_2^{\theta_2} \left( \int_0^1 \left( \int_0^1 \left( t_1^{\frac{1}{p_1}+\varepsilon_1} t_2^{\frac{1}{p_2}+\varepsilon_2} f^{*_1,*_2}(t_1,t_2) \right)^{q_1} \frac{dt_1}{t_1} \right)^{\frac{q_2}{q_1}} \frac{dt_2}{t_2} \right)^{\frac{1}{q_2}},$$

және $0 < \bar{\delta} < \frac{1}{\bar{p}}$ үшін:

$$\|f\|_{GL_{\bar{p},\bar{q}}^{-\bar{\theta}}(\Omega)} \asymp \inf_{\substack{0 < \varepsilon_1 \le \delta_1 \\ 0 < \varepsilon_2 \le \delta_2}} \varepsilon_1^{\theta_1} \varepsilon_2^{\theta_2} \left( \int_0^1 \left( \int_0^1 \left( t_1^{\frac{1}{p_1}-\varepsilon_1} t_2^{\frac{1}{p_2}-\varepsilon_2} f^{*_1,*_2}(t_1,t_2) \right)^{q_1} \frac{dt_1}{t_1} \right)^{\frac{q_2}{q_1}} \frac{dt_2}{t_2} \right)^{\frac{1}{q_2}}.$$

**Дәлелдеу.** Теореманы дәлелдеу үшін $\delta_1, \delta_2 < 1, \eta_1, \eta_2 > 1, \eta_1 \delta_1 = 1, \eta_2 \delta_2 = 1, \sigma_1 = \varepsilon_1 \eta_1, \sigma_2 = \varepsilon_2 \eta_2$ екенін ескереміз:

$$\sup_{\substack{0 < \varepsilon_1 \le \delta_1 \\ 0 < \varepsilon_2 \le \delta_2}} \varepsilon_1^{\theta_1} \varepsilon_2^{\theta_2} \left( \int_0^1 \left( \int_0^1 \left( t_1^{\frac{1}{p_1}+\varepsilon_1} t_2^{\frac{1}{p_2}+\varepsilon_2} f^{*_1,*_2}(t_1,t_2) \right)^{q_1} \frac{dt_1}{t_1} \right)^{\frac{q_2}{q_1}} \frac{dt_2}{t_2} \right)^{\frac{1}{q_2}}$$

$$= \sup_{\substack{0 < \frac{\sigma_1}{\eta_1} \le \delta_1 \\ 0 < \frac{\sigma_2}{\eta_2} \le \delta_2}} \left(\frac{\sigma_1}{\eta_1}\right)^{\theta_1} \left(\frac{\sigma_2}{\eta_2}\right)^{\theta_2} \left( \int_0^1 \left( \int_0^1 \left( t_1^{\frac{1}{p_1}+\frac{\sigma_1}{\eta_1}} t_2^{\frac{1}{p_2}+\frac{\sigma_2}{\eta_2}} f^{*_1,*_2}(t_1,t_2) \right)^{q_1} \frac{dt_1}{t_1} \right)^{\frac{q_2}{q_1}} \frac{dt_2}{t_2} \right)^{\frac{1}{q_2}}$$

$$\le \eta_1^{-\theta_1} \eta_2^{-\theta_2} \sup_{\substack{0 < \sigma_1 \le \eta_1 \delta_1 \\ 0 < \sigma_2 \le \eta_2 \delta_2}} \sigma_1^{\theta_1} \sigma_2^{\theta_2} \left( \int_0^1 \left( \int_0^1 \left( t_1^{\frac{1}{p_1}+\frac{\sigma_1}{\eta_1}} t_2^{\frac{1}{p_2}+\frac{\sigma_2}{\eta_2}} f^{*_1,*_2}(t_1,t_2) \right)^{q_1} \frac{dt_1}{t_1} \right)^{\frac{q_2}{q_1}} \frac{dt_2}{t_2} \right)^{\frac{1}{q_2}}$$

$$\leq \eta_1^{-\theta_1}\eta_2^{-\theta_2} \sup_{\substack{0<\sigma_1\leq 1 \\ 0<\sigma_2\leq 1}} \sigma_1^{\theta_1}\sigma_2^{\theta_2} \left(\int_0^1 \left(\int_0^1 \left(t_1^{\frac{1}{p_1}+\sigma_1} t_2^{\frac{1}{p_2}+\sigma_2} f^{*_1,*_2}(t_1,t_2)\right)^{q_1} \frac{dt_1}{t_1}\right)^{\frac{q_2}{q_1}} \frac{dt_2}{t_2}\right)^{\frac{1}{q_2}}$$

$$= C\|f\|_{GL_{\bar{p},\bar{q}}^{\bar{\theta}}(\Omega)}$$

Айталық $0 < \delta_1 < \frac{1}{p_1}$, $0 < \delta_2 < \frac{1}{p_2}$, $\eta_1, \eta_2 > 1$, $\eta_1\delta_1 = \frac{1}{p_1}$, $\eta_2\delta_2 = \frac{1}{p_2}$, $\sigma_1 = \varepsilon_1\eta_1$, $\sigma_2 = \varepsilon_2\eta_2$ болсын:

$$\inf_{\substack{0<\varepsilon_1\leq \delta_1 \\ 0<\varepsilon_2\leq \delta_2}} \varepsilon_1^{\theta_1}\varepsilon_2^{\theta_2} \left(\int_0^1 \left(\int_0^1 \left(t_1^{\frac{1}{p_1}-\varepsilon_1} t_2^{\frac{1}{p_2}-\varepsilon_2} f^{*_1,*_2}(t_1,t_2)\right)^{q_1} \frac{dt_1}{t_1}\right)^{\frac{q_2}{q_1}} \frac{dt_2}{t_2}\right)^{\frac{1}{q_2}}$$

$$= \inf_{\substack{0<\frac{\sigma_1}{\eta_1}\leq \delta_1 \\ 0<\frac{\sigma_2}{\eta_2}\leq \delta_2}} \left(\frac{\sigma_1}{\eta_1}\right)^{\theta_1}\left(\frac{\sigma_2}{\eta_2}\right)^{\theta_2} \left(\int_0^1 \left(\int_0^1 \left(t_1^{\frac{1}{p_1}-\frac{\sigma_1}{\eta_1}} t_2^{\frac{1}{p_2}-\frac{\sigma_2}{\eta_2}} f^{*_1,*_2}(t_1,t_2)\right)^{q_1} \frac{dt_1}{t_1}\right)^{\frac{q_2}{q_1}} \frac{dt_2}{t_2}\right)^{\frac{1}{q_2}}$$

$$\leq \eta_1^{-\theta_1}\eta_2^{-\theta_2} \inf_{\substack{0<\sigma_1\leq \frac{1}{p_1} \\ 0<\sigma_2\leq \frac{1}{p_2}}} \sigma_1^{\theta_1}\sigma_2^{\theta_2} \left(\int_0^1 \left(\int_0^1 \left(t_1^{\frac{1}{p_1}-\sigma_1} t_2^{\frac{1}{p_2}-\sigma_2} f^{*_1,*_2}(t_1,t_2)\right)^{q_1} \frac{dt_1}{t_1}\right)^{\frac{q_2}{q_1}} \frac{dt_2}{t_2}\right)^{\frac{1}{q_2}}$$

$$= C\|f\|_{GL_{\bar{p},\bar{q}}^{-\bar{\theta}}}.$$

**Теорема 6.** Егер $\bar{\theta} < \bar{\lambda}$, $\bar{\lambda}_1 - \bar{\theta}_1 = \frac{1}{q_1} - \frac{1}{\tau_1}$, $\bar{\theta}_2 = \frac{1}{q_2} - \frac{1}{\tau_2}$, $0 < p_1, p_2 < \infty$ болса, онда келесі енгізу орынды:

$$GL_{\bar{p},\bar{\tau}}^{\bar{\theta}}(\Omega) \hookrightarrow GL_{\bar{p},\bar{q}}^{\bar{\lambda}}(\Omega).$$

**Дәлелдеу.** Теореманы дәлелдеу үшін Гельдер теңсіздігін қолданамыз:

$$\|f\|_{GL_{\bar{p},\bar{q}}^{\bar{\lambda}}(\Omega)} = \sup_{0<\varepsilon_1,\varepsilon_2<1} \varepsilon_1^{\lambda_1}\varepsilon_2^{\lambda_2} \left( \int_0^1 \left( \int_0^1 \left( t_1^{\frac{1}{p_1}+\varepsilon_1} t_2^{\frac{1}{p_2}+\varepsilon_2} f^{*_1,*_2}(t_1,t_2) \right)^{q_1} \frac{dt_1}{t_1} \right)^{\frac{q_2}{q_1}} \frac{dt_2}{t_2} \right)^{\frac{1}{q_2}}$$

$$\leq \sup_{0<\varepsilon_1,\varepsilon_2<1} \varepsilon_1^{\lambda_1}\varepsilon_2^{\lambda_2} \left( \int_0^1 \left( \int_0^1 \left( t_1^{\frac{1}{p_1}+\frac{\varepsilon_1}{2}} t_2^{\frac{1}{p_2}+\frac{\varepsilon_2}{2}} f^{*_1,*_2}(t_1,t_2) \right)^{q_1 h} \frac{dt_1}{t_1} \right)^{\frac{q_2}{q_1}} \frac{dt_2}{t_2} \right)^{\frac{1}{q_2 h}}$$

$$\cdot \left( \int_0^1 \left( \int_0^1 \left( t_1^{\frac{\varepsilon_1}{2}} t_2^{\frac{\varepsilon_2}{2}} \right)^{q_1 h'} \frac{dt_1}{t_1} \right)^{\frac{q_2}{q_1}} \frac{dt_2}{t_2} \right)^{\frac{1}{q_2 h'}}.$$

Мұндағы $q_1 h = \tau_1$ және $q_2 h = \tau_2$ деп таңдаймыз. Теңсіздіктің екінші көбейтіндісін қарастырайық:

$$\left( \int_0^1 \left( \int_0^1 \left( t_1^{\frac{\varepsilon_1}{2}} t_2^{\frac{\varepsilon_2}{2}} \right)^{q_1 h'} \frac{dt_1}{t_1} \right)^{\frac{q_2}{q_1}} \frac{dt_2}{t_2} \right)^{\frac{1}{q_2 h'}} = \frac{1}{\left|\frac{\varepsilon_1}{2} q_1 h'\right|^{\frac{1}{q_1 h'}} \cdot \left|\frac{\varepsilon_2}{2} q_2 h'\right|^{\frac{1}{q_2 h'}}} = C \frac{1}{\varepsilon_1^{\frac{1}{q_1 h'}} \cdot \varepsilon_2^{\frac{1}{q_2 h'}}}.$$

Яғни:

$$\|f\|_{GL_{\bar{p},\bar{q}}^{\bar{\lambda}}(\Omega)} \leq \sup_{0<\varepsilon_1,\varepsilon_2<1} \varepsilon_1^{\lambda_1}\varepsilon_2^{\lambda_2} \left( \int_0^1 \left( \int_0^1 \left( t_1^{\frac{1}{p_1}+\frac{\varepsilon_1}{2}} t_2^{\frac{1}{p_2}+\frac{\varepsilon_2}{2}} f^{*_1,*_2}(t_1,t_2) \right)^{q_1 h} \frac{dt_1}{t_1} \right)^{\frac{q_2}{q_1}} \frac{dt_2}{t_2} \right)^{\frac{1}{q_2 h}}$$

$$\cdot C \frac{1}{\varepsilon_1^{\frac{1}{q_1 h'}} \cdot \varepsilon_2^{\frac{1}{q_2 h'}}}$$

$$\leq C \sup_{0<\varepsilon_1,\varepsilon_2<1} \varepsilon_1^{\lambda_1 - \frac{1}{q_1 h'}} \varepsilon_2^{\lambda_2 - \frac{1}{q_2 h'}} \left( \int_0^1 \left( \int_0^1 \left( t_1^{\frac{1}{p_1}+\frac{\varepsilon_1}{2}} t_2^{\frac{1}{p_2}+\frac{\varepsilon_2}{2}} f^{*_1,*_2}(t_1,t_2) \right)^{q_1 h} \frac{dt_1}{t_1} \right)^{\frac{q_2 h}{q_1 h}} \frac{dt_2}{t_2} \right)^{\frac{1}{q_2 h}}$$

$$\leq C \sup_{0<\varepsilon_1,\varepsilon_2<1} \varepsilon_1^{\theta_1}\varepsilon_2^{\theta_2} \left( \int_0^1 \left( \int_0^1 \left( t_1^{\frac{1}{p_1}+\varepsilon_1} t_2^{\frac{1}{p_2}+\varepsilon_2} f^{*,*_2}(t_1,t_2) \right)^{\tau_1} \frac{dt_1}{t_1} \right)^{\frac{\tau_2}{\tau_1}} \right)^{\frac{1}{\tau_2}} = C\|f\|_{GL_{\bar{p},\bar{\tau}}^{\bar{\theta}}(\Omega)}.$$

**Теорема 7.** Егер $0 < \bar{p} < \infty, 0 < \bar{\tau} \leqslant \infty$ және $\bar{\theta} > 0$ болса, онда

$$\|f\|_{GL_{\bar{p},\bar{\tau}}^{\bar{\theta}}(\Omega)} \asymp$$

$$\sup_{1<k_1,k_2} k_1^{-\theta_1} k_2^{-\theta_2} \left( \sum_{m_1=-\infty}^{0} \left( \sum_{m_2=-\infty}^{0} \left( 2^{m_1\left(\frac{1}{p_1}+\frac{1}{k_1}\right)} 2^{m_2\left(\frac{1}{p_2}+\frac{1}{k_2}\right)} f^{*_1,*_2}(2^{m_1}, 2^{m_2}) \right)^{\tau_1} \right)^{\frac{\tau_2}{\tau_1}} \right)^{\frac{1}{\tau_2}}.$$

**Дәлелдеу.**

$$\|f\|_{GL_{\bar{p},\bar{\tau}}^{\bar{\theta}}(\Omega)} = \sup_{0<\varepsilon_1,\varepsilon_2<1} \varepsilon_1^{\theta_1} \varepsilon_2^{\theta_2} \left( \int_0^1 \left( \int_0^1 \left( t_1^{\frac{1}{p_1}+\varepsilon_1} t_2^{\frac{1}{p_2}+\varepsilon_2} f^{*_1,*_2}(t_1,t_2) \right)^{\tau_1} \frac{dt_1}{t_1} \right)^{\frac{\tau_2}{\tau_1}} \frac{dt_2}{t_2} \right)^{\frac{1}{\tau_2}}.$$

$$\varepsilon_1^{\theta_1} \varepsilon_2^{\theta_2} \left( \int_0^1 \left( \int_0^1 \left( t_1^{\frac{1}{p_1}+\varepsilon_1} t_2^{\frac{1}{p_2}+\varepsilon_2} f^{*_1,*_2}(t_1,t_2) \right)^{\tau_1} \frac{dt_1}{t_1} \right)^{\frac{\tau_2}{\tau_1}} \frac{dt_2}{t_2} \right)^{\frac{1}{\tau_2}} =$$

$$= \left( \sum_{k_1=-\infty}^{0} \left( \sum_{k_2=-\infty}^{0} \int_{2^{k_1-1}}^{2^{k_1}} \int_{2^{k_2-1}}^{2^{k_2}} \left( t_1^{\frac{1}{p_1}+\varepsilon_1} t_2^{\frac{1}{p_2}+\varepsilon_2} f^{*_1,*_2}(t_1,t_2) \right)^{\tau_1} \frac{dt_1}{t_1} \right)^{\frac{\tau_2}{\tau_1}} \frac{dt_2}{t_2} \right)^{\frac{1}{\tau_2}} =$$

$$= \left( \sum_{k_1=0}^{\infty} \left( \sum_{k_2=0}^{\infty} \int_{2^{-k_1-1}}^{2^{-k_1}} \int_{2^{-k_2-1}}^{2^{-k_2}} \left( t_1^{\frac{1}{p_1}+\varepsilon_1} t_2^{\frac{1}{p_2}+\varepsilon_2} f^{*_1,*_2}(t_1,t_2) \right)^{\tau_1} \frac{dt_1}{t_1} \right)^{\frac{\tau_2}{\tau_1}} \frac{dt_2}{t_2} \right)^{\frac{1}{\tau_2}}.$$

$$\left( \int_{2^{-k_1-1}}^{2^{-k_1}} \left( \int_{2^{-k_2-1}}^{2^{-k_2}} \left( t_1^{\frac{1}{p_1}+\varepsilon_1} t_2^{\frac{1}{p_2}+\varepsilon_2} f^{*_1,*_2}(t_1,t_2) \right)^{\tau_1} \frac{dt_1}{t_1} \right)^{\frac{\tau_2}{\tau_1}} \frac{dt_2}{t_2} \right)^{\frac{1}{\tau_2}}$$

$$\begin{cases} \leq (f^{*_1,*_2}(2^{-k_1-1}, 2^{-k_2-1}) 2^{-k_1\left(\frac{1}{p_1}+\varepsilon_1\right)} 2^{-k_2\left(\frac{1}{p_2}+\varepsilon_2\right)})^{\tau_1} \\ \geq (f^{*_1,*_2}(2^{-k_1}, 2^{-k_2}) 2^{(-k_1-1)\left(\frac{1}{p_1}+\varepsilon_1\right)} 2^{(-k_2-1)\left(\frac{1}{p_2}+\varepsilon_2\right)})^{\tau_2}. \end{cases}$$

$$I \leq \left( \sum_{k_1=0}^{\infty} \left( \sum_{k_2=0}^{\infty} \left( \left( 2^{-k_1\left(\frac{1}{p_1}+\varepsilon_1\right)} 2^{-k_2\left(\frac{1}{p_2}+\varepsilon_2\right)} \right) f^{*_1,*_2}(2^{-k_1-1}, 2^{-k_2-1})^{\tau_1} \right)^{\frac{\tau_2}{\tau_1}} \right)^{\frac{1}{\tau_2}},$$

$$I \geq \left( \sum_{k_1=0}^{\infty} \left( \sum_{k_2=0}^{\infty} \left( f^{*_1,*_2}(2^{-k_1}, 2^{-k_2}) 2^{(-k_1-1)\left(\frac{1}{p_1}+\varepsilon_1\right)} 2^{(-k_2-1)\left(\frac{1}{p_2}+\varepsilon_2\right)} \right)^{\tau_1} \right)^{\frac{\tau_2}{\tau_1}} \right)^{\frac{1}{\tau_2}}.$$

$\|f\|_{GL_{\vec{p},\vec{\tau}}^{-\vec{\theta}}(\Omega)}$ жағдайында супремум инфинумге ауысады:

$$\|f\|_{GL_{\vec{p},\vec{\tau}}^{-\vec{\theta}}(\Omega)} \asymp \inf_{1<k_1,k_2} k_1^{\theta_1} k_2^{\theta_2} \left( \sum_{m_1=-\infty}^{0} \left( \sum_{m_2=-\infty}^{0} \left( 2^{m_1\left(\frac{1}{p_1}-\frac{1}{k_1}\right)} 2^{m_2\left(\frac{1}{p_2}-\frac{1}{k_2}\right)} f^{*_1,*_2}(2^{m_1}, 2^{m_2}) \right)^{\tau_1} \right)^{\frac{\tau_2}{\tau_1}} \right)^{\frac{1}{\tau_2}}$$

**Қорытынды**

Бұл мақалада жаңа анизотропты гранд Лоренц кеңістіктері анықталып, олардың қасиеттері зерттелді. Бұл кеңістіктер функционалдық кеңістіктерді зерттеу үшін біртұтас параметрлік құрылымды ұсынады, әсіресе параметрлердің шектік жағдайларын зерттеуге мүмкіндік береді. Мақалада анизотропты гранд Лоренц кеңістіктерінің анизотропты Лоренц кеңістіктерінен айырмашылықтары көрсетіліп, олардың енгізу теоремалары, нормаланған кеңістіктер және функцияларды зерттеу мәселелеріне қатысты қасиеттері талқыланды.

Мақалада келесі негізгі нәтижелер алынды:
1. Гранд анизотропты Лоренц кеңістіктерінің анықтамасы;
2. Бағалаулар;
3. Енгізу теоремалары;
4. Параметрлер бойынша енгізу шарттары.

Бұл зерттеу Гранд анизотропты Лоренц кеңістіктерінің қасиеттерін тереңірек түсінуге мүмкіндік береді. Осы кеңістіктердің параметрлік шкаласы әртүрлі функционалдық кеңістіктерді зерттеуге және математикалық есептерді шешуге арналған жаңа әдіс-тәсілдерді әзірлеуге негіз болуы мүмкін. Болашақта бұл кеңістіктердің қолданылуы функционалдық талдау саласындағы зерттеулерге үлкен үлес қосады.

**Әдебиеттер**

[1,2]**Manarbek M.**
PhD Student, ORCID ID: 0009-0006-6879-8356,
e-mail: makpal9136@mail.ru
[2]**Tleukhanova N. T.**
Professor, Doctor of Physical and Mathematical Sciences,
ORCID ID: 0000-0002-4133-7780,
e-mail: tleukhanova63@rambler.ru
[2]**Mussabayeva G. K.**
PhD, ORCID ID: 0000-0003-2368-8955,
e-mail: mussabayeva@mail.ru

Institute of Mathematics and Mathematical Modeling, Almaty, Kazakhstan,
L.N. Gumilyov Eurasian National University, Astana, Kazakhstan.

## ANISOTROPIC GRAND LORENTZ SPACES AND THEIR PROPERTIES

**Abstract**

In this article, new anisotropic grand Lorentz spaces are defined and their properties are studied. These spaces are a new structure that provides a unified parameter for the study of various functional spaces. The consideration of grand spaces is especially important for the study of boundary conditions of parameters and allows us to achieve new results in this area. The study of boundary parameters in classical spaces is not always possible. In recent years, grand Lebesgue spaces and their generalizations have been widely studied in problems of functional spaces. These spaces are generalizations of classical Lorentz and grand Lorentz spaces. The article defines grand anisotropic Lorentz spaces, gives basic estimates in these spaces, proves embedding theorems, and derives embedding theorems for parameters. The results obtained can play an important role not only in theoretical, but also in applied problems.

**Key words:** Lorentz spaces, grand Lorentz spaces, embedding theorems, inequalities, anisotropic spaces.

[1,2]**Манарбек М.**,
PhD студент, ORCID ID: 0009-0006-6879-8356,
e-mail: manarbek@math.kz

[2]**Тлеуханова Н. Т.**,
профессор, доктор физико-математических наук,
ORCID ID: 0000-0002-4133-7780,
e-mail: tleukhanova@rambler.ru

[2]**Мусабаева Г. К.**,
PhD, ORCID ID: 0000-0003-2368-8955,
e-mail: musabaevaguliya@mail.ru

[1]Институт математики и математического моделирования, Алматы, Казахстан,
[2]Евразийский национальный университет им. Л.Н. Гумилева, Астана, Казахстан

# АНИЗОТРОПНЫЕ ГРАНД ПРОСТРАНСТВА ЛОРЕНЦА И ИХ СВОЙСТВА

**Аннотация**

В данной статье определяются новые анизотропные гранд пространства Лоренца и изучаются их свойства. Эти пространства представляют собой новую структуру, которая обеспечивает единую среду для исследования различных функциональных пространств. Рассмотрение гранд пространств особенно важно для изучения граничных условий параметров, и в этом отношении могут быть получены новые результаты. Не всегда возможно изучить граничные параметры в классических пространствах. В последние годы гранд пространства Лебега и их обобщения широко изучаются в задачах функциональных пространств. Эти пространства являются обобщениями классических пространств Лоренца и больших пространств Лоренца. В статье определяются большие анизотропные пространства Лоренца, приводятся основные оценки в этих пространствах, доказываются теоремы вложения и выводятся теоремы вложения для параметров. Полученные результаты могут сыграть важную роль не только в теоретическом плане, но и в прикладных задачах.

**Ключевые слова:** пространства Лоренца, гранд пространства Лоренца, теоремы вложения, неравенства, анизотропные пространства.